\documentclass[12pt]{amsart}
\usepackage{amsmath,amssymb,amsfonts}

\newtheorem{OP}{Open Problem}
\theoremstyle{plain}
\newtheorem{thm}{Theorem}[section]

\newtheorem{conj}[thm]{Conjecture}

\theoremstyle{definition}
\newtheorem{df}[thm]{Definition}
\newtheorem{rem}[thm]{Remark}

\newcommand\bR{{\mathbb R}}
\newcommand\bE{{\mathbb E}}
\newcommand\bP{{\mathbb P}}
\newcommand\bZ{{\mathbb Z}}
\newcommand\bC{{\mathbb C}}

\newcommand\Sing{\operatorname{Sing}}

\newcommand\gt{\theta}

\newcommand\T{\Theta}
\newcommand\e{\varepsilon}
\newcommand\de{\delta}

\newcommand\cA{{\mathcal A}}
\newcommand\cM{{\mathcal M}}
\newcommand\cJ{{\mathcal J}}
\newcommand\cH{{\mathcal H}}

\newcommand\cP{{\mathcal P}}
\newcommand\cS{{\mathcal S}}
\newcommand\cU{{\mathcal U}}
\newcommand\cI{{\mathcal I}}

\newcommand\Pic{\operatorname{Pic}}
\newcommand\Sp{{\operatorname{Sp}(2g,\bZ)}}

\newcommand\tc[2]{{\gt\left[\begin{matrix}{#1}\\{#2}\end{matrix}\right]}}
\begin{document}
\title{The Schottky Problem}
\author{Samuel Grushevsky}
\address{Mathematics Department, Stony Brook University,
Stony Brook, NY 11790-3651, USA}
\email{sam@math.sunysb.edu}
\thanks{Research supported in part by National Science Foundation under the grant DMS-0901086/DMS-1053313.}

\begin{abstract}
In this survey we discuss some of the classical and modern methods in studying the (Riemann-)Schottky problem, the problem of characterizing Jacobians of curves among principally polarized abelian varieties. We present many of the recent results in this subject, and describe some directions of current research. This paper is based on the talk given at the ``Classical algebraic geometry today'' workshop at MSRI in January 2009.
\end{abstract}

\maketitle

\section{Introduction}
The Riemann-Schottky problem is the problem of determining which complex principally polarized abelian varieties arise as Jacobian varieties of complex curves. The history of the problem is very long, going back to the works of Abel, Jacobi, and Riemann. The first approach, culminating in a complete solution in genus 4 (the first non-trivial case), was developed by Schottky \cite{schottky} and Schottky-Jung \cite{scju}. Since then a variety of different approaches to the problem have been developed, and many geometric properties of abelian varieties in general and Jacobian varieties in particular have been studied extensively. Numerous partial and complete solutions to the Schottky problem have been conjectured, and some were proven.

%However, a solution in the classical spirit of Schottky is still not known even for genus 5, and for a given 5-dimensional abelian variety there is no known way to definitively determine whether it is a Jacobian or not.

In this survey we will describe many of the ideas and methods that have been applied to or developed for the study of the Schottky problem. We will present some of the results, as well as various open problems and possible connections among various approaches. To keep the length of the text reasonable, the proofs for the most part will be omitted, and references will be given; when possible, we will try to indicate the general idea or philosophy behind the work done. We hope that an interested reader may consider this as an introduction to the ideas and results of the subject, and would be encouraged to explore the field in greater depth by following some of the references.

This text is in no way the first (and will certainly not be the last) survey written on the Schottky problem. Many excellent surveys, from various points of view, and emphasizing various aspects of the field, have been written, including those by Dubrovin \cite{dusurvey}, Donagi \cite{dosurvey}, Beauville \cite{besurvey}, Debarre \cite{desurvey}, Taimanov \cite{tasurvey}, van Geemen \cite{vgsurvey}, Arbarello \cite[appendix]{mumfarb}, Buchstaber and Krichever \cite{bukrsurvey}. A beautiful introduction is \cite{mumfordbookjacobians}, while many relevant results on curves, abelian varieties, and theta functions can be found in \cite{igusabook,mumfordbooktheta1,mumfordbooktheta2,mumfordbooktheta3,acgh,bila}.  The research literature on or related to the Schottky problem is vast, with many exciting results dating back to the works of Schottky, and more progress constantly happening. While we have strived to give as many references as possible, our list is by no means complete, and we apologize for any inadvertent omissions; another good source of references is in Arbarello's appendix to \cite{mumfarb}

\smallskip
The structure of this work is as follows. In section \ref{intro} we introduce most of the notation for curves, abelian varieties, and their moduli, and state the problem; in section \ref{classical} we introduce theta constants and discuss the classical approach culminating in an explicit solution of the Schottky problem in genus 4; in section \ref{schottkyjung} we continue the discussion of modular forms vanishing on the Jacobian locus, and describe the Schottky-Jung approach to obtaining them using the Prym construction. In section \ref{andreottimayer} we discuss the singularities of the theta divisor for Jacobians and Pryms, and present the Andreotti-Mayer approach to the Schottky problem. In section \ref{minimalclass} we discuss minimal cohomology classes and the Matsusaka-Ran criterion; in section \ref{kummerimage} we discuss the geometry of the Kummer variety, its secants, and in particular the trisecant conjecture characterizing Jacobians, and the Prym analog. Section \ref{gamma00} deals with the $\Gamma_{00}$ conjecture and further geometry of the $|2\T|$ linear system.

\subsection*{Acknowledgements}
This paper grew out of a talk given at the ``Classical algebraic geometry today'' workshop during the Algebraic Geometry program at MSRI in January 2009. I am grateful to the organizers of that conference for the opportunity to present the talk and especially to Mircea Musta\c t\u a and Mihnea Popa for the invitation and encouragement to write this survey. I would also like to thank Igor Krichever, Mihnea Popa, and Riccardo Salvati Manni for many conversations over the years, from which I learned a lot about different viewpoints and results on the Schottky problem. I am grateful to Olivier Debarre, Mihnea Popa, and the referee for carefully reading the manuscript and for numerous useful suggestions and advice. We thank Enrico Arbarello and Yuri Zarhin for drawing our attention to their works in \cite{mumfarb},\cite{zarhin1},\cite{zarhin2}.

\section{Notation: the statement of the Schottky problem}\label{intro}
Loosely speaking, the Schottky problem is the following question: which principally polarized abelian varieties are Jacobians of curves? In this section we introduce the relevant notation and state this question more precisely. Such a brief introduction cannot do justice to the interesting and deep constructions of the moduli stacks of curves and abelian varieties; we will thus assume that the reader is in fact familiar with the ideas of moduli theory, take the existence of moduli space for granted, and will use this section primarily to fix notation.

Throughout the text we work over the field $\bC$ of complex numbers. Unless stated otherwise, we are working in the category of smooth projective varieties.

\begin{df}[Complex tori and abelian varieties]
A $g$-dimensional complex torus is a quotient of $\bC^g$ by a full rank lattice, i.e.~by a subgroup $\Lambda\subset\bC^g$ (under addition) such that $\Lambda\sim\bZ^{2g}$ and $\Lambda\otimes_{\bZ}\bR=\bR^{2g}=\bC^g$.

A $g$-dimensional complex torus is called an abelian variety if it is a projective variety. This is equivalent to the torus admitting a very ample line bundle, and it turns out that this is equivalent to the lattice $\Lambda$ being conjugate under the action of $GL(g,\bC)$ to a lattice $\bZ^g+\tau\bZ^g$ for some symmetric complex $g\times g$ matrix $\tau$ with positive-definite imaginary part. This statement is known as Riemann's bilinear relations, see for example \cite{bila}. Note that for a given $\Lambda$ there exist many possible $\tau$ such that $\Lambda$ is isomorphic to $\bZ^g+\tau\bZ^g$ --- this is made more precise below.
\end{df}

\begin{df}[Polarization]
A polarization on an abelian variety $A$ is the first Chern class of an ample line bundle $L$ on $A$, i.e. the polarization is $[L]:=c_1(L)\in H^2(A,\bZ)\cap H^{1,1}(A,\bC)$. A polarization is called principal if $h^0(A,L)=1$ (in terms of the cohomology class, this means that the matrix is unimodular). We will denote principally polarized abelian varieties (abbreviated as ppavs) by $(A,\Theta)$, where $\Theta$ is our notation for the ample bundle, the first Chern class of which is the polarization.
\end{df}

We note that the polarizing line bundle $L$ on an abelian variety can be translated by an arbitrary point of the abelian variety, preserving the class of the polarization $[L]$. However, it is customary to use the bundle rather than its class in the notations for a ppav. Given a polarization, one usually chooses a symmetric polarizing bundle, i.e~one of the $2^{2g}$ line bundles $\Theta$ (differing by a translation by a 2-torsion point) of the given Chern class $[\Theta]$ such that $(-1)^*\Theta=\Theta$ --- where $-1$ is the involution of the abelian variety as a group.

\begin{df}[Siegel space]
We denote by $\cH_g$ the space of all symmetric $g\times g$ matrices $\tau$ with positive definite imaginary part --- called the Siegel space. Given an element $\tau\in\cH_g$, the Riemann theta function
is defined by the following Fourier series
\begin{equation}\label{thetadef}
  \theta(\tau,z):=\sum\limits_{n\in\bZ^g} \exp(\pi i n^t\tau n+2\pi in^t z).
\end{equation}
\end{df}

This (universal) theta function is a map $\theta:\cH^g\times\bC^g\to\bC$. For fixed $\tau\in\cH^g$ it has the following easily verified from the definition automorphy property in $z$:
\begin{equation}\label{ztrans}
 \theta(\tau,z+\tau n+ m)=\exp(-\pi i n^t\tau n-2\pi i n^t z)\theta(\tau,z)\quad\forall n,m\in\bZ^g
\end{equation}
It thus follows that the zero locus $\lbrace z\in\bC^g\mid \theta (\tau,z)=0\rbrace$ is invariant under the shifts by the lattice $\bZ^g+\tau\bZ^g$, and thus descends to a well-defined subvariety $\Theta_\tau\subset\bC^g/(\bZ^g+\tau\bZ^g$).

\begin{df}[Moduli of ppavs]
We denote by $\cA_g$ the moduli stack of principally polarized (complex) abelian varieties $(A,\Theta)$ of dimension $g$. Then we have a natural map $\cH_g\to\cA_g$ given by
$$
  \tau\mapsto \left(A_\tau:=\bC^g/(\bZ^g+\tau\bZ^g,)\Theta_\tau:=\lbrace\theta(\tau,z)=0\rbrace\right).
$$
It can be shown that this map exhibits $\cH_g$ as the universal cover of $\cA_g$. Moreover, there is a natural action of the symplectic group $\Sp$ on $\cH_g$ given by
\begin{equation}\label{Spaction}
 \gamma\circ \tau:=(\tau c+d)^{-1}(\tau a+b),
\end{equation}
where we write $\Sp\ni\gamma=\begin{pmatrix} a&b\\ c&d\end{pmatrix}$ in $g\times g$ block form. It turns out that we have $\cA_g=\cH_g/\Sp$, so that the stack $\cA_g$ is in fact a global quotient of a smooth complex manifold by a discrete group. For future use we denote
$$
 \cA_g^{dec}:=\left(\mathop{\cup}\limits_{i=1\ldots g-1}\cA_i\times\cA_{g-i}\right)\subset\cA_g
$$
the locus of decomposable ppavs (products of lower-dimensional ones) --- these turn out to be special in many ways, as these are the ones for which the theta divisor is reducible.

We also denote by $\cA_g^{ind}:=\cA_g\setminus \cA_g^{dec}$ its complement, the locus of indecomposable ppavs .
We refer for example to \cite{bila} for more details on $\cA_g$, and to \cite{igusabook,mumfordbooktheta1,mumfordbooktheta2,mumfordbooktheta1} for the theory of theta functions.
\end{df}
\begin{df}[Universal family]
We denote by $\cU_g\to\cA_g$ the universal family of ppavs, the fiber over $(A,\T)$ being the ppav itself. Analytically, the universal cover of $\cU_g$ is $\cH_g\times\bC^g$, since the universal cover of any ppav is $\bC^g$. We thus have
$$
 \cU_g=\cH_g\times\bC^g/(\Sp\ltimes\bZ^{2g}),
$$
where $\bZ^{2g}$ acts by fixing $\cH_g$ and adding lattice vectors on $\bC^g$, while $\gamma\in\Sp$ acts on $\bC^g$ by $z\mapsto (\tau c+d)^{-1}z$.
\end{df}
\smallskip
\begin{df}[Jacobians]
Analytically, to define the Jacobian of a Riemann surface $C$ of genus $g$, one considers some standard generators $A_1,\ldots, A_g$, $B_1,\ldots, B_g$ of the fundamental group $\pi_1(C)$ (i.e.~generators such that the only relation is $\prod_{i=1}^gA_iB_iA_i^{-1}B_i^{-1}=1$). Then one chooses a basis $\omega_1,\ldots,\omega_g$ of the space of abelian differentials $H^0(C,K_C)$ dual to $\lbrace A_i\rbrace$, i.e.~such that $\int_{A_i}\omega_j=\delta_{i,j}$. The Jacobian is then the ppav $J(C):=(A_\tau,\T_\tau)$ corresponding to the matrix $\lbrace\tau_{ij}\rbrace:=\lbrace \int_{B_i}\omega_j\rbrace\in\cH_g$.

Algebraically, the Jacobian of a curve $C$ of genus $g$ is $J(C):=\Pic^{g-1}(C)$ --- the set of linear equivalence classes of divisors of degree $g-1$ on $C$. The principal polarization divisor on it is the set of effective bundles $\Theta:=\lbrace L\in\Pic^{g-1}(C)\mid h^0(C,L)\ge 1\rbrace$.
\end{df}

The Jacobian is (non-canonically) isomorphic to $\Pic^n(C)$; an isomorphism is defined by choosing a fixed divisor of degree $g-1-n$ and subtracting it. Under this identification $\Pic^n(C)$ has a natural principal polarization, but, unlike $\Pic^{g-1}(C)$, there is no natural choice of a polarization divisor. It is often convenient to view the Jacobian as $\Pic^0(C)$, which is naturally an abelian group (while $\Pic^{g-1}(C)$ is naturally a torsor over $\Pic^0(C)$ rather than a group in a natural way).

\begin{df}[Moduli of curves]
We denote by $\cM_g$ the moduli stack of smooth (compact complex) curves of genus $g$, and by $\overline{\cM_g}$ its Deligne-Mumford compactification \cite{demu} --- the moduli stack of stable curves of genus $g$. The Torelli map is the map $J:\cM_g\to\cA_g$ sending $C$ to its Jacobian $J(C)$. The Torelli theorem states that this map is injective, and we call its image the Jacobian locus  $\cJ_g$ (it is sometimes called the Torelli locus in the literature).
\end{df}

We recall that $\dim\cM_g=3g-3$ (for $g>1$), while $\dim\cA_g=g(g+1)/2$. Thus the dimensions coincide for $g\le 3$, and in fact the Jacobian locus $\cJ_g$ is equal to $\cA_g^{ind}$ iff $g\le 3$. For $g=4$ we have $\dim\cM_4=9=10-1=\dim\cA_4-1$, and $\cJ_4$ has been completely described by Schottky \cite{schottky} (his result was proven rigorously by Igusa \cite{igusagen4} and Freitag \cite{freitaggen4}, and further reformulated by Igusa \cite{igusachristoffel}).

The subject of this survey is the following problem, studied by Riemann in unpublished notes, and with first published significant progress made by Schottky \cite{schottky} and Schottky-Jung \cite{scju}:

\smallskip
\noindent {\bf (Riemann)-Schottky problem:}\\
Describe the image $\cJ_g:=J(\cM_g)\subset\cA_g$.

\smallskip
In many of the approaches to the Schottky problem that we will survey weak solutions to this problem will be obtained --- i.e.~we will get loci within $\cA_g$, of which $\cJ_g$ would be an irreducible component.

Note that the moduli space $\cM_g$ is not compact, and neither is the Jacobian locus $\cJ_g$. However, it was shown by Namikawa \cite{namikawa1,namikawa2} that the Torelli map extends to a map $\overline J:\overline{\cM_g}\to\overline{\cA_g}$ from the Deligne-Mumford compactification of $\cM_g$ to a certain (second Voronoi, but we will not need the definition here) toroidal compactification of $\cA_g$. Thus in some sense it is more natural to consider the closed Schottky problem: that of describing the closure $\overline J_g\subset\overline{\cA_g}$.

A lot of work has been done on the Schottky problem, and much progress has been made. Outwardly it may even appear that the problem has been solved completely in many ways. While we survey the many solutions or partial solutions that have indeed been developed, it will become apparent that a lot of work is still to be done, and much substantial understanding is still lacking, with many important open questions remaining. In particular, the following two questions are completely open:

\begin{OP}\label{SchottkyOP} Solve the Schottky problem explicitly for $g=5$, i.e.~given an explicit $\tau\in\cH_5$, determine whether $[\tau]\in\cA_5$ lies in $\cJ_5$.
\end{OP}

A lot of research has been done on $\cA_5$, and many partial solutions to the Schottky problem for $g=5$ are available, see \cite{accola,donagiscju,casalaina2,cmfr}, but a complete solution in the spirit of Schottky (similar to genus 4, see the next section) is still not known; there is no known way to answer the above question. Another question, which would be a potential application for a solution to the Schottky problem, is:
\begin{OP}\label{ColemanOP}
The symmetric space $\cH_g$ has a natural $\operatorname{Sp}(2g,\bR)$-invariant metric on it, which descends to a metric on $\cA_g$. For $g$ sufficiently large do there exist any complex geodesics for this metric contained in $\overline{\cJ_g}$ and intersecting $\cJ_g$?
\end{OP}
A negative answer to this question, together with the Andr\'e-Oort conjecture \cite{andrebook,oortconj}, would imply the Coleman conjecture \cite{coleman}, stating that for $g$ sufficiently large
there exist only finitely many $C\in\cM_g$ such that $J(C)$ admits a complex multiplication.

If we had a good explicit solution to the Schottky problem it could allow one to check explicitly whether such geodesics may exist, but so far this has not been accomplished, though much interesting work in the area has been done \cite{toledo,hain,vizu,MVZ}.

\smallskip
Questions related to the Schottky problem also seem to arise naturally in perturbative (super)string theory, and many deep mathematical questions were posed and studied by D'Hoker-Phong \cite{dhp1,dhp2,dhp3,dhp4,dhpazyg1,dhpazyg2}, with many more still open questions arising in recent works \cite{cdpvg1,gstring,smstring, cdpvg2,grsmgen5,opsy} and others.

The Schottky problem as we stated it is still a loosely phrased question --- there may be different ways of describing $\cJ_g$. In the following sections we survey the many approaches to characterizing Jacobians.

\section{Theta constants: the classical approach}\label{classical}
The classical approach to the Schottky problem, initiated in the works of Riemann and Schottky, is to try to embed $\cA_g$ in a projective space, and then try to write out the defining ideal for the image of $\cJ_g$. If successful, this would of course be a very explicit solution to the Schottky problem.

\begin{df}[Hodge bundle and modular forms]
The Hodge rank $g$ vector bundle $\bE$ over $\cA_g$ has the space of holomorphic one-forms $H^0(A,\Omega^1_A)=H^{1,0}(A)$ as the fiber over a ppav $A$. Its determinant line bundle $L:=\det\bE$ is  called the bundle of (scalar) modular forms of weight 1, and in general we call a section of $L^{\otimes k}$ a modular form of weight $k$.

The sections of any line bundle on $\cA_g$ can be pulled back to $\cH_g$ (where the pullback of the bundle is trivial, as $\cH_g$ is contractible). Thus the sections of any line bundle on $\cA_g$ can be alternatively described as holomorphic functions on $\cH_g$ satisfying certain automorphy properties under the action of the group of deck transformation $\Sp$. In particular, we have
\begin{equation}\label{modform}
 H^0(\cA_g,L^{\otimes k})=\left\lbrace f:\cH_g\to\bC\mid f(\gamma\circ\tau)=\det(\tau c+d)^k f(\tau)\right\rbrace
\end{equation}
for all $\gamma\in\Sp$, written in $g\times g$ block form, and for all $\tau\in\cH_g$, with the action given by (\ref{Spaction}). This is the analytic definition of modular forms, by their automorphy properties.
\end{df}
It can be shown (this is work of Igusa and Mumford, see \cite{bila} for more details) that $L$ is an ample line bundle on $\cA_g$, and thus its sufficiently high power gives an embedding of $\cA_g$ into a projective space. However, the actual very ample power is quite high, and instead it is easier to work with level covers of $\cA_g$, i.e.~with certain finite covers of $\cA_g$.
\begin{df}[Level covers]
Algebraically, for any $\ell\in\bZ_{>0}$ we denote by $\cA_g(\ell)$ the full level $\ell$ cover of $\cA_g$, i.e.~the moduli space of ppavs $(A,\Theta)$ together with a chosen symplectic (with respect to the Weil pairing) basis for the group $A[\ell]$ of points of order $\ell$ on the abelian variety. We denote by $\cA_g(\ell,2\ell)$ the theta level $\ell$ cover of $\cA_g$, i.e.~the moduli space of ppavs together with a chosen symplectic basis for $A[\ell]$ and also one chosen point of order $2\ell$. For $\ell$ divisible by two both of these are Galois covers, with the Galois groups described as follows.

The principal level $\ell$ subgroup of $\Sp$ is defined as
$$
  \Gamma_g(\ell):=\left\lbrace\gamma=\begin{pmatrix} a&b\\
  c&d\end{pmatrix} \in\Sp \right|\left.\, \gamma\equiv \begin{pmatrix} 1&0 \\   0&1\end{pmatrix}\ {\rm mod}\ \ell\right\rbrace
$$
while the theta level $\ell$ subgroup is defined as
$$
 \Gamma_g(\ell,2\ell):=\left\lbrace \gamma\in\Gamma_g(\ell)\mid
  {\rm diag} (a^tb)\equiv{\rm diag}(c^td)\equiv0\ {\rm mod}\ 2\ell\right\rbrace.
$$
These are the subgroups preserving the level data as above, i.e.~we have $\cA_g(\ell)=\cH_g/\Gamma_g(\ell)$ and $\cA_g(\ell,2\ell)=\cH_g/\Gamma_g(\ell,2\ell)$.

One then defines the level Jacobian loci $\cJ_g(\ell)$ and $\cJ_g(\ell,2\ell)$ as the preimages of $\cJ_g$ under the corresponding level covers of $\cA_g$.
\end{df}
The usefulness of these covers lies in the fact that it is much easier to construct modular forms on them (i.e.~sections of the pullbacks of $L$ to them). In particular the following construction gives a large supply of modular forms of weight $1/2$.
\begin{df}[Theta functions with characteristics]
For $\e,\de\in \frac{1}{\ell}\bZ^g/\bZ^g$ (which we will also think of as a point $m=\tau\e+\de\in A_\tau[\ell]$) the
theta function with characteristic $\e,\de$ (or $m$) is
\begin{equation}\label{thetachar}
\begin{aligned}
\tc\e\de(\tau,z):=\theta_m(\tau,z):&=
\sum\limits_{n\in\bZ^g} \exp\left(\pi i(n+\e)^t\left(\tau(n+\e)+2(z+\de)\right)\right)\\ &=\exp\left(\pi i\e^t\tau\e+2\pi i\e^t(z+\de)\right)\theta(\tau,z+m)
\end{aligned}
\end{equation}
As a function of $z$,  $\theta_m(\tau,z)$ is a section of the translate of the principal polarization $\T_\tau$ by the point $m$.

Setting $z=0$, we get theta constants of order $\ell$, denoted $\theta_m(\tau):=\theta_m(\tau,0)$, which by above formula are, up to some easy exponential factor, just the values of the Riemann theta function at points $m$ of order $\ell$.
\end{df}

It turns out that theta constants $\theta_m(\tau)$ are modular forms of weight $1/2$ for  $\Gamma_g(2\ell,4\ell)$, i.e.~are sections in $H^0(\cA_g(2\ell,4\ell),L^{\otimes(1/2)})$ (the pullback to this level cover of the bundle $L$ of modular forms has a square root, i.e.~a bundle such that its
tensor square is $L$ --- this is what we mean by $L^{\otimes(1/2)}$). Moreover, they define an embedding of the moduli space:
\begin{thm}[Igusa \cite{igusabook} for $\ell=4n^2$, Mumford for $\ell\ge 3$,
Salvati Manni \cite{smlevel2} for $\ell\ge 2$]
The map $\Phi_\ell:\cA_g(2\ell,4\ell)\to\bP^{\ell^{2g}-1}$ given by $
\tau\mapsto\left\lbrace \theta_m(\tau)\right\rbrace_{\forall m\in A_\tau[\ell]}$ is an embedding. Thus the bundle $L^{\otimes(1/2)}$ is very ample on $\cA_g(2\ell,4\ell)$ (In fact theta constants with characteristics generate the space of its sections, and generate the ring of modular forms on $\cA_g(2\ell,4\ell)$).
\end{thm}

From now on we will concentrate on the simplest case, that of level $\ell=2$. In this case the theta function $\theta_m(\tau,z)$ is an even (resp.~odd) function of $z$ if $(-1)^{4\de^t\e}=1$ (resp.~$=-1)$; we will call such $m\in A[2]$ even (resp.~odd). Thus the the constant $\theta_m(\tau)$ vanishes identically if (and in fact only if) $m\in A[2]^{odd}$. There are $2^{g-1}(2^g+1)$ even theta constants, and thus we can consider the image $\Phi_2(\tau)\in\bP^{2^{g-1}(2^g+1)-1}$; however, this is still a large space, and even for $g=1$ the map $\Phi_2$ is not dominant --- the one defining equation for the image is in fact $\tc00(\tau)^4=\tc01(\tau)^4+\tc10(\tau)^4$. For $g=2$ the situation already gets very complicated  (see \cite{vgeemenscju,vgvdg} for more details). Thus it is natural to try to embed a suitable cover of $\cA_g$ in a smaller projective space.

\begin{df}[Theta functions of the second order]
For $\e\in\frac12\bZ^g/\bZ^g$ the theta function of the second order with characteristics $\e$ is
\begin{equation}\label{2Tdef}
  \T[\e](\tau,z):=\theta_{2\tau\e}(2\tau,2z)
\end{equation}
\end{df}
For $\tau$ fixed we have $\T[\e](\tau,z)\in H^0(A_\tau,2\T_\tau)$, and it can be shown that in fact theta functions of the second order generate the space of sections of $2\T_\tau$. Noting that for any $m=\tau\e+\de\in A_\tau[2]$ the square of the theta function with characteristic $\theta_m^2(\tau,z)$ for $\tau$ fixed is also a section of $2\T_\tau$, we can express it as a linear combination of $\T[e](\tau,z)$. The result is the following identity
\begin{thm}[Riemann's bilinear addition formula]
\begin{equation}\label{bilinearidentity}
 \theta_m^2(\tau,z)=\sum\limits_{\sigma\in \frac12\bZ^g/\bZ^g} (-1)^{4\delta^t\sigma} \T[\sigma+\epsilon](\tau,0)\T[\sigma](\tau,z).
\end{equation}
\end{thm}
\begin{df}[The theta map]
We define the theta constants of the second order to be $\T[\e](\tau):=\T[\e](\tau,0)$; it turns out that these are also modular forms of weight $1/2$, but for $\Gamma_g(2,4)$, i.e.~$\T[\e](\tau)\in H^0(\cA_g(2,4),L^{\otimes (1/2)})$.
We can thus consider the theta map
\begin{equation}\label{Th}
 Th:\tau\to\lbrace \T[\e](\tau)\rbrace_{\forall \e\in\frac12\bZ^g/\bZ^g}.
\end{equation}
\end{df}
Notice that Riemann's bilinear addition formula above, for $z=0$, shows that $\Phi_2(\tau)$ can be recovered uniquely up to signs from $Th(\tau)$. Since $\Phi_2$ is injective on $\cA_g(4,8)$, it follows that the map $Th$ is finite-to-one on $\cA_g(4,8)$, and thus also finite-to-one on $\cA_g(2,4)$, its natural domain of definition. In fact it is known that $Th:\cA_g(2,4)\to\bP^{2^g-1}$ is generically injective, and conjecturally it is an embedding \cite{smlevel2}. Thus the classical formulation of the Schottky problem is the following question:

\smallskip
\noindent {\bf Classical Riemann-Schottky problem:}\\
Write the defining equations for $\overline{Th(\cJ_g(2,4))}\subset \overline{Th(\cA_g(2,4))}\subset\bP^{2^g-1}.$

\smallskip
Notice that for $g=1,2$ the target projective space is of the same dimension as the moduli, and it is in fact known that $Th$ is then dominant. For $g=3$ the 6-dimensional image $\overline{Th(\cA_g(2,4))}=\overline{Th(\cJ_3(2,4))}$ is a hypersurface in $\bP^7$. The defining equation for it can be written easily in terms of polynomials of even degree in theta constants with characteristics --- which can then be rewritten in terms of $Th$ by using Riemann's bilinear addition formula. Miraculously the same kind of formula gives the one defining equation for $\overline{Th(\cJ_4(2,4))}\subset \overline{Th(\cA_4(2,4))}$:

\begin{thm}[Schottky \cite{schottky}, Igusa \cite{igusachristoffel}, Farkas-Rauch \cite{farascju}, see also van Geemen-van der Geer \cite{vgvdg}]
In genus $g$ define the Igusa modular form to be
\begin{equation}
 F_g(\tau):=2^g\sum\limits_{m\in A[2]}\theta_m^{16}(\tau) -\left(\sum\limits_{m\in A[2]}\theta_m^{8}(\tau)\right)^2.
\end{equation}
This is a modular form with respect to the entire group $\Sp$, such that when rewritten in terms of theta constants of the second order using Riemann's bilinear addition formula (\ref{bilinearidentity}), the form $F_g$ is
\begin{itemize}
\item identically zero for $g=1,2$;
\item the defining equation for  $\overline{Th(\cJ_3(2,4))}=\overline{Th(\cA_3(2,4))}\subset\bP^7$ for $g=3$;
\item the defining equation for $\overline{Th(\cJ_4(2,4))}\subset\overline{Th(\cA_4(2,4))}$ for $g=4$.
\end{itemize}
\end{thm}
To date this remains the most explicit solution to the Schottky problem in genus 4, and no similar solution in higher genera is known or has been proposed. However, the following question, which could be a test case for tackling open problem \ref{ColemanOP} has not been settled:
\begin{OP}
Construct all geodesics for the invariant metric on $\cA_4$ contained in $\overline{\cJ_4}$, and intersecting $\cJ_4$.
\end{OP}

\smallskip
One of course wonders what happens for $g>4$. It turns out that the form $F_g$ also naturally arises in perturbative superstring theory \cite{dhpazyg1,grsmgen5}. By interpreting $F_g$ in terms of lattice theta functions, $F_g=\theta_{D_{16}^+}-\theta_{E_8\times E_8}$ (see \cite{conwaysloanebook} for definitions), using some physics intuition, and with motivation from the works of Belavin, Knizhnik, Morozov, it was conjectured by D'Hoker and Phong \cite{dhpazyg1} that $F_g$ vanishes on $\cJ_g$ for any genus. If true, this would have been a very nice defining equation for the Jacobian locus (though of course more would be needed). However, the situation is more complicated.

\begin{thm}[Grushevsky-Salvati Manni \cite{grsmgen5}]
The modular form $F_g$ does not vanish identically on $\cJ_g$ for any $g\ge 5$. In fact the zero locus of $F_5$ on $\cJ_5$ is the locus of trigonal curves.
\end{thm}
It was shown by Poor \cite{poorFg} that for any $g$ the form $F_g$ vanishes on the locus of Jacobians of hyperelliptic curves, and in view of the above it is natural to ask whether it also vanishes on the trigonal locus. As a result, the following easier step on the way to solving open problem \ref{SchottkyOP} still has not been accomplished:

\begin{OP}
Write at least one nice/invariant modular form vanishing on $\cJ_5$.
\end{OP}
\begin{rem}
Many modular forms vanishing on $\cJ_5$ were produced by Accola \cite{accola} by using the Schottky-Jung approach (see the next section of this survey). Accola showed that these modular forms give a weak solution to the Schottky problem, i.e.~define a locus in $\overline{Th(\cA_5(2,4))}$ of which the Jacobian locus is an irreducible component. However, these modular forms are not ``nice'' in the sense that a lot of combinatorics of theta characteristics is involved in obtaining them, and they are certainly modular forms for $\Gamma_5(2,4)$, not for all of $\operatorname{Sp}(10,\bZ)$.
\end{rem}

\begin{rem}
The codimension of $\cJ_5$ within $\cA_5$ is equal to $3$, and it can be shown (Faber \cite{faberalgorithms}) that $\overline {Th(\cJ_5(2,4))}\subset\overline {Th(\cA_5(2,4))}$ is not a complete intersection. In fact the degrees of these subvarieties of $\bP^{31}$ can be computed \cite{grdeg}. One would then expect that it may be possible to check whether Accola's equations produce extra components besides the Jacobian locus, but since the defining ideal of $\overline {Th(\cA_5(2,4))}\subset\bP^{31}$ is not known, at the moment this seems impossible.
\end{rem}

Instead of asking to characterize the Jacobian locus, one could ask different Schottky-type questions --- to characterize the image in $\cA_g$ of certain special subvarieties of $\cM_g$. In sharp contrast to the original Schottky problem, the hyperelliptic Schottky problem admits a complete solution, via simple explicit equations for theta constants:
\begin{thm}[Mumford \cite{mumfordbooktheta2}, Poor \cite{poor}]
For any $g$ there exist some (explicitly described in terms of certain combinatorics) sets of characteristics $S_1,\ldots, S_N\subset \frac12\bZ^{2g}/\bZ^{2g}$ such that $\tau\in\cA_g$ is the period matrix of a hyperelliptic Jacobian if and only if for some $1\le i\le N$ we have
$$
\theta_m(\tau)=0\Longleftrightarrow m\in S_i\qquad \forall m\in\frac12\bZ^{2g}/\bZ^{2g}.
$$
We remark that this solution even takes care of the locus of decomposable ppavs (for them more vanishing happens).
\end{thm}

\section{Modular forms vanishing on the Jacobian locus: the Schottky-Jung approach}\label{schottkyjung}
Despite the fact that at the moment a solution to the classical Schottky problem in genus 5 is not known, there is a general method to produce many modular forms vanishing on the Jacobian locus. This is the original approach originated by Schottky-Jung \cite{scju}, and developed by Farkas-Rauch \cite{farascju}.

\begin{df}[Prym variety]
For an \'etale connected double cover $\widetilde{C}\to C$ of a curve $C\in\cM_g$ --- such a cover is given by a two-torsion point $\eta\in J(C)[2]\setminus\lbrace0\rbrace$ --- we define the Prym variety to be
$$
 Prym(\widetilde C\to C):=Prym(C,\eta):=Ker_0(J(\widetilde{C})\to J(C))\in\cA_{g-1},
$$
where $Ker_0$ denotes the connected component of 0 in the kernel (which in fact has two connected components). The map $J(\widetilde C)\to J(C)$ here is the norm map corresponding to the cover $\widetilde C\to C$, so the Prym is the connected component of zero in the kernel of the norm map.

The restriction of the principal polarization $\T_{J(\widetilde C)}$ to the Prym gives twice the principal polarization. However, it turns out that this polarization admits a canonical square root, which thus gives a natural principal polarization on the Prym. We denote by $\cP_{g-1}\subset\cA_{g-1}$ the locus of Pryms of all \'etale covers of curves in $\cM_g$; the question of describing it is called the Prym-Schottky problem. See Mumford \cite{mumfordprym} for a modern exposition of the Prym construction.
\end{df}
The idea of the Schottky-Jung approach is as follows: the curve $\widetilde{C}$, being an \'etale double cover, is special (i.e.~ for $g>1$ the set of such double covers is a proper subvariety of $\cM_{2g-1}$). It can be shown that this implies the vanishing of some theta constants for $J(\widetilde C)$, see \cite{farascju}. On the other hand, there exists a finite-to-one surjective map (such maps are called isogenies of abelian varieties) $\left(J(C)\times Prym(\widetilde C\to C)\right)\longrightarrow J(\widetilde C)$, and pulling back by it allows one to express the theta function on $J(\widetilde C)$ in terms of theta functions on $J(C)$ and on the Prym. Combining this with the vanishing properties for theta constants on $J(\widetilde C)$ gives the following:

\begin{thm}[Schottky-Jung proportionality \cite{scju}, proven rigorously by Farkas-Rauch \cite{farkashscju}]
Let $\tau$ be the period matrix of $C$ and let $\pi$ be the period matrix of the Prym for $\eta=\left[\begin{matrix}0&0&\ldots&0\\
1&0&\ldots&0\end{matrix}\right]$. Then for any $\e,\de\in\frac12 \bZ^{g-1}/\bZ^{g-1}$ the theta constants of $J(C)$ and of the Prym are related by
\begin{equation}\label{scju}
\tc\e\de(\pi)^2= c\tc{0\,\e}{0\,\de}(\tau)\cdot\tc{0\,\e}{1\,\de}(\tau)
\end{equation}
where $\pi$ is the period matrix of the Prym, $\tau$ is the period matrix of $J(C)$, and the constant $c$ is independent of $\e,\de$.
\end{thm}
To obtain the versions of Schottky-Jung proportionalities for double covers corresponding to other points $\eta$ we can act on the above equation by $\Sp/\Gamma_g(1,2)$, which acts transitively on the set of two-torsion points, and permutes the theta constants up to certain eighth roots of unity (see \cite{igusabook,bila} for precise theta transformation formulas). Analytically the situation for $\eta=\left[\begin{matrix}0&0&0&\ldots&0\\
1&1&0&\ldots&0\end{matrix}\right]$ was worked out by Farkas \cite{farkashscju} --- an extra sign depending on $\e,\de$ is present then.

\begin{df}[The Schottky-Jung locus]
Let $I_{g-1}$ be the defining ideal for the image $\overline{Th(\cA_{g-1}(2,4))}\subset\bP^{2^{g-1}-1}$. For any equation $F\in I_{g-1}$ let $F_\eta$ be the polynomial equation on $\bP^{2^g-1}$ obtained by using the Schottky-Jung proportionality to substitute an appropriate polynomial of degree 2 in terms of theta constants of $\tau$ for the square of any theta constant of $\pi$ (as a result square roots of polynomials in theta constants may appear, and we would then take the product of the expressions for all possible choices of the square root, i.e.~the norm).

Let then $S_g^\eta$ be the ideal obtained from $I_{g-1}$ in this way. The (big) Schottky-Jung locus $\cS_g^\eta(2,4)\subset\cA_g(2,4)$ is then defined to be the zero locus of $S_g^\eta$. Note that it is a priori not clear --- and in fact not known --- that $I_g\subset S_g^\eta$, and thus it may make a difference that we define $\cS_g^\eta(2,4)$ within $\cA_g(2,4)$, and not as a subvariety of the projective space. We then define the small Schottky-Jung locus to be
\begin{equation}\label{schottkylocus}
\cS_g(2,4):=\mathop{\bigcap}\limits_{\eta\in\frac12\bZ^{2g}/\bZ^{2g}\setminus\lbrace 0\rbrace}\cS_g^\eta(2,4).
\end{equation}
Since the action of $\Sp$ permutes the different $\eta$ and the ideals $S_g^\eta$, it follows that the ideal defining $\cS_g(2,4)$ is $\Sp$-invariant, and thus the locus $\cS_g(2,4)$ is a preimage of some $\cS_g\subset\cA_g$ under the level cover.
\end{df}

\begin{thm}[van Geemen \cite{vgeemenscju}, Donagi \cite{donagiscju}]
The Jacobian locus $\cJ_g$ is an irreducible component of the small Schottky-Jung locus $\cS_g$, and in fact $\cJ_g(2,4)$ is an irreducible component of the big Schottky-Jung locus $\cS_g^\eta(2,4)$ for any $\eta$.
\end{thm}
In genus 4 it can in fact be shown that $\cS_4^\eta(2,4)=\cJ_4(2,4)$ for one (and thus for all) $\eta$ (see \cite{igusagen4}), and thus also $\cS_4=\cJ_4$. However, in genus 5 it turns out that $\cJ_5(2,4) \subsetneq\cS_5^\eta(2,4)$, since the latter contains the locus of intermediate Jacobians --- see \cite{donagiintermjac}. However, it can be shown that the locus of intermediate Jacobians does not lie in $\cS_5$, and the following bold conjecture could be made.
\begin{conj}[Donagi \cite{dosurvey}]\label{sjconj}
The small Schottky-Jung locus is equal to the Jacobian locus, $\cS_g=\cJ_g$.
\end{conj}

\begin{rem}
The dimension of the space of \'etale double covers of curves of genus $g+1$ is equal to $3g$, the same as $\dim\cM_{g+1}$, and thus $\dim\cP_g\le 3g$. It follows that for $g\ge 6$ the Prym locus $\cP_g$ is a proper subvariety of $\cA_g$ (it is in fact known that $\cP_g$ is dense within $\cA_g$ for $g\le 5$, see \cite{dosm}). Thus for $g\ge 7$ one could apply the Schottky-Jung proportionality to the defining ideal for $\overline{Th(\cP_{g-1})}\subsetneq\overline{ Th(\cA_{g-1})}$ to get further equations for $\overline{Th(\cJ_g)}$. Thus the conjecture above would seem more natural for $g=5$ and perhaps $g=6$ than in general. On the other hand, the approaches in \cite{vgeemenscju,donagiscju} to the Schottky-Jung theory have been via degeneration, and in that case questions related to the $\Gamma_{00}$ conjecture (see section \ref{kummerimage}) appear on the boundary; it could be that the recent progress on that topic could shed more light on the Schottky-Jung locus.
\end{rem}

\smallskip
To summarize, the classical Schottky-Jung approach to the Schottky results in the explicit  solution for $g=4$, a ``weak'' solution (i.e.~up to extra components) in any genus, and conjecturally could yield a complete solution in every genus --- assuming the ideal of $\overline{Th(\cA_g(2,4))}$ were known completely (many elements of this ideal are known, but we do not know that they generate the entire ideal, see \cite{vgeemenscju,smrelations}).

\section{Singularities of the theta divisor: the Andreotti-Mayer approach}\label{andreottimayer}
It was shown by Andreotti and Mayer \cite{anma} that for a generic ppav the theta divisor is smooth. The singularities of the theta divisor for a Jacobian of a curve can be described using
\begin{thm}[Riemann's theta singularity theorem; see \cite{acgh}]
For any curve $C$ and any $D\in J(C)=\Pic^{g-1}(C)$ we have ${\rm mult}_D\T=h^0(C,D)$.
\end{thm}
We note that this theorem agrees with the algebraic definition of the theta divisor $\T$ on the Jacobian as the locus of effective divisors). This theorem was further generalized by Kempf \cite{kempf}. From this theorem we see that $\Sing\T_{J(C)}=\lbrace D\in\Pic^{g-1}\mid h^0(C,D)\ge 2\rbrace$. It follows that the theta divisor for Jacobians of curves has a large singular set:
\begin{thm}[Andreotti-Mayer \cite{anma}]
For a non-hyperelliptic curve $C$ of genus $g$ we have $\dim(\Sing\T_{J(C)})=g-4$, while for a hyperelliptic curve $C$ this dimension is $g-3$.
\end{thm}

One would thus expect that this property is special for Jacobians. However, one could not expect that it defines $\overline{\cJ_g}\subset\overline{\cA_g}$. Indeed, note that for a decomposable ppav $(A,\T)=(A_1,\T_1)\times(A_2,\T_2)$ (where $(A_i,\T_i)\in\cA_{g_i}$ with $g_1+g_2=g$) we have $\T=(\T_1\times A_2)\cup (A_1\times \T_2)$, and thus $\Sing\T\supset \T_1\times \T_2$ is of dimension $g-2$. It is thus natural to study the following
\begin{df}[Andreotti-Mayer loci]
We define the $k$'th Andreotti-Mayer locus to be
$$N_{k,g}:=\lbrace (A,\T)\in\cA_g \mid \dim\Sing\T\ge k\rbrace.$$
\end{df}
By definition we have $N_{k,g}\subset N_{k+1,g}$. Of course we have $N_{g-1,g}=\emptyset$, and by the above we see that  $\cA_g^{dec}\subset N_{g-2,g}$. It was conjectured by Arbarello-De Concini \cite{ardcnovikov} and proven by Ein-Lazarsfeld \cite{eila1} that in fact $N_{g-2,g}=\cA_g^{dec}$. Thus one is led to study the next cases, which include the Jacobians:
\begin{thm}[Andreotti-Mayer \cite{anma}]
$\cJ_g$ is an irreducible component of $N_{g-4,g}$, while the locus of hyperelliptic Jacobians $Hyp_g$ is an irreducible component of $N_{g-3,g}$.
\end{thm}
Thus we obtain a weak solution to the Schottky problem, and it is natural to ask to describe $(N_{g-4,g}\cap\cA_g^{ind})\setminus \cJ_g$ and $(N_{g-3,g}\cap\cA_g^{ind})\setminus Hyp_g$. In the simplest cases we have $N_{0,3}\cap\cA_3^{ind}=Hyp_3$, while Beauville \cite{beauville} showed that $N_{0,4}$ in fact has two irreducible components, one of which is $\cJ_4$, and the other one is the theta-null divisor:
$$
 \theta_{{\rm null},g}:=\Big\lbrace\tau\ \mid\prod\limits_{m\in A[2]^{even}}\theta_m(\tau)=0\Big\rbrace=\left\lbrace (A,\T)\in\cA_g \mid A[2]^{\rm even}\cap\T\ne\emptyset\right\rbrace.
$$

Since in the universal family of ppavs $\cU_g\to\cA_g$ the locus of (fiberwise) singularities of the theta divisor is given by $g+1$ equations $\theta(\tau,z)=\frac\partial{\partial z_1}\theta(\tau,z)=\ldots=\frac\partial{\partial z_g}\theta(\tau,z)=0$, this locus has expected codimension $g+1$ in $\cU_g$. Thus $N_{0,g}$,  its image under $\cU_g\to\cA_g$ (with fiber dimension $g$), has expected codimension 1. It was shown by Beauville \cite{beauville} that $N_{0,g}$ is indeed always a divisor, and it was shown by Debarre \cite{debarredecomposes} that for any $g\ge 4$ this divisor has two irreducible components, scheme-theoretically $N_{0,g}=\theta_{{\rm null},g}+2N'_{0,g}$.

One possible complete solution to the Schottky problem would be the following
\begin{conj}[see Debarre \cite{debarrecodim3}]
The following set-theoretic equalities hold:
$$N_{g-3,g}\cap\cA_g^{ind}=Hyp_g;\qquad N_{g-4,g}\setminus\cJ_g\subset\theta_{{\rm null},g}.$$
\end{conj}
In view of this conjecture it is natural to try to understand the intersection of $\overline{\cJ_g}$ with the other irreducible components of $N_{g-4,g}$ and, to start with, to study $\overline{\cJ_g}\cap\theta_{{\rm null},g}$.
%By studying the universal family of singularities of theta divisors over $\overline{\cJ_g}$ one can deduce certain special properties for ppavs in this intersection.
The simplest interesting case is that of $N_{0,4}=\theta_{{\rm null},4}\cup\overline{\cJ_4}$. The following result, conjectured by H.~Farkas \cite{farkashvanishing}, was proven
\begin{thm}[Grushevsky-Salvati Manni\cite{grsmgen4}, Smith-Varley\cite{smvagen4}]
The locus of Jacobians of curves of genus 4 with a vanishing theta-null is equal to the locus of 4-dimensional ppavs for which the double point singularity of the theta divisor is not ordinary (i.e.~the tangent cone does not have maximal rank): this is to say
$$\begin{aligned}\overline{\cJ_4}\cap\gt_{{\rm null},4}&=\lbrace A\in\cA_4\mid \exists m\in A[2]^{\rm even};\gt(\tau,m)=\det_{i,j}\frac{\partial^2\gt(\tau,z)}{\partial z_i\partial z_j}|_{z=m}=0\rbrace\\
&=\lbrace A\in\cA_4\mid\exists m\in A[2]^{\rm even}\cap\T;
\ TC_m\T {\rm\ has\ rank\ } \le 3\rbrace \end{aligned}$$
\end{thm}
It is thus natural to denote the locus above by $\gt_{{\rm null},4}^{3}$, for rank of the tangent cone being at most 3. In general one has
\begin{thm}[Debarre \cite{debarredecomposes}, Grushevsky-Salvati Manni \cite{grsmordertwo}, Smith-Varley \cite{smvagen4}]
$$(\overline{\cJ_g}\cap\gt_{{\rm null},g})\ \subset\ \gt_{{\rm null},g}^3\ \subset\ \gt_{{\rm null},g}^{g-1}\ \subset\ (\gt_{{\rm null},g}\cap N_0')\subset\Sing N_0$$
\end{thm}
Here the first inclusion follows from Kempf's \cite{kempf} generalization of Riemann theta singularity theorem, the second inclusion is obvious by definition, the third inclusion is the content of the theorem, and the fourth inclusion is immediate since the singular locus of a reducible variety contains the intersections of its irreducible components.

At the moment this theorem seems to be the best result in trying to understand the conjecture above, and many further questions remain.

\begin{rem}
The Andreotti-Mayer loci are also of importance for the Prym-Schottky problem. There is an analog of Riemann's theta singularity theorem for Prym varieties (see \cite{casalaina}), showing that their theta divisor are also very singular, and the following result holds:
\begin{thm}[Debarre \cite{debarreprymam}]
The Prym locus $\cP_g$ is an irreducible component of $N_{g-6,g}$.
\end{thm}
It is thus also natural to try to understand the other components of $N_{g-6,g}$. Debarre (\cite{debarrecodim3}, corollary 12.5) constructed some irreducible components of this locus, but even for the simplest case of $g=6$ the locus $N_{0,6}$ is not completely understood.
\end{rem}

\smallskip
One of the main difficulties in following this approach to the Schottky problem is the difficulty in understanding the Andreotti-Mayer loci $N_{k,g}$. In fact even their dimensions are not known, though there is the following
\begin{conj}[Ciliberto-van der Geer \cite{civdg1,civdg2}]
For any $1\le k\le g-3$ and for any irreducible component $X\subset N_{k,g}$ such that for $(A,\T)\in X$ general $\mathop{End}(A,\T)=\bZ$ (in particular $X\not\subset\cA_g^{dec}$), we have
$${\rm codim}_{\cA_g} X\ge \frac{(k+1)(k+2)}{2},$$
with equality only for components $\cJ_g\subset N_{g-4,g}$ and $Hyp_g\subset N_{g-3,g}$.
\end{conj}
Unfortunately still not much is known about this conjecture: it was shown by Mumford \cite{mumforddimag} that $N_{1,g}\subsetneq N_{0,g}$, so that ${\rm codim} N_{1,g}\ge 2$. Ciliberto and van der Geer's \cite{civdg1,civdg2} proof of this conjecture for $k=1$ (which shows that ${\rm codim} N_{1,g}\ge 3$) is already very hard. To the best of our knowledge even the following question is open:
\begin{OP}
Can it happen that $N_{k,g}=N_{k+1,g}$ for some $k,g$?
\end{OP}

\begin{rem}
We would like to mention here a related question, which does not seem to shed any new light on the Schottky problem, but is of interest in the study of $\cA_g$. Indeed, one can study other properties of $\T$; and in particular instead of looking at $\dim\Sing\T$ one could define the multiplicity loci
$$
 S_k:=\lbrace (A,\T)\in\cA_g\mid \exists z\in A, {\rm mult}_z\T\ge k\rbrace
$$
It was shown by Koll\'ar \cite{kollarbook} that $S_{g+1}=\emptyset$, while Smith-Varley \cite{smva} showed that $S_g$ is equal to the locus of products of elliptic curves $S_g=\cA_1\times\ldots\times\cA_1$ .

However, from Riemann's theta singularity theorem it follows that the maximum multiplicity of the theta function on Jacobians is $\lfloor \frac{g+1}{2}\rfloor$, and it can be shown that this is also the case for Pryms (see \cite{casalaina}), i.e.~that
$$
  \cJ_g\cap S_{\lfloor \frac{g+3}{2}\rfloor}=\cP_g\cap S_{\lfloor \frac{g+3}{2}\rfloor}=\emptyset.
$$
Since one expects the theta divisors of Jacobians and Prym in general to be very singular, this leads to the following natural
\begin{conj}
$$\cA_g^{ind}\cap S_{\lfloor \frac{g+3}{2}\rfloor}=\emptyset.$$
\end{conj}
For $g\le 5$ this conjecture was shown by Casalaina-Martin \cite{casalaina} to be true, using the fact that then the Prym locus $\cP_g$ is dense in $\cA_g$. We refer to \cite{civdg2,grsmconjectures} for more details and further conjectures on the loci of ppavs with different sorts of singularities of the theta divisor.

To summarize, the Andreotti-Mayer approach gives geometric conditions for a ppav to be a Jacobian, and the original results of Andreotti and Mayer give a geometric weak solution to the Schottky problem (though in practice it is not easy to compute $\dim\Sing\T_\tau$ for an explicitly given $\tau\in\cH_g$). With recent progress, the Andreotti-Mayer approach can be improved to give a complete solution to the Schottky problem in genus 4, since it gives a geometric description of $N_{0,4}\setminus\overline{\cJ_4}$ as $\theta_{{\rm null},4}\setminus\theta_{{\rm null},4}^3$. However, already for the $g=5$ case the situation is not completely understood, though there are compelling conjectures for arbitrary genus. The relationship between the Andreotti-Mayer approach and other properties of the theta divisor were studied in much more detail by Beauville-Debarre \cite{bede} and Debarre \cite{debarrestratification}.

\section{Subvarieties of a ppav: minimal cohomology classes}\label{minimalclass}
One striking approach to the Schottky problem stems from the observation that for a Jacobian $J(C)$ one can naturally map the symmetric product $Sym^d C$ (for $1\le d< g$) to $J(C)=\Pic^{g-1}(C)$ by fixing a divisor $D\in\Pic^{g-1-d}(C)$ and mapping $(p_1,\ldots, p_d)\mapsto D+\sum p_i$. The image of such a map, denoted $W^d(C)\subset J(C)$, is independent of $D$ up to a translation, and one can compute its cohomology class
$$
 [W^d(C)]=\frac{[\T]^{g-d}}{(g-d)!}\in H^{2g-2d}(J(C)),
$$
where by $[\T]$ we denote the cohomology class of the polarization.
One can show that this cohomology class is indivisible in cohomology with $\bZ$ coefficients, and we thus call this class minimal. Of course these subvarieties of $J(C)$ are very special --- in particular $W^1(C)\simeq C$, and one can thus ask whether their existence is a special property of Jacobians. There is a following criterion
\begin{thm}[Matsusaka \cite{matsusaka}, Ran \cite{ranmin1}]
A ppav $(A,\T)$ is a Jacobian if and only if there exists a curve $C\subset A$ with $[C]=\frac{[\T]^{g-1}}{(g-1)!}$ (in which case $(A,\T)=J(C)$).
\end{thm}
This gives a complete geometric solution to a weaker form of the Schottky problem: given a pair $C\subset A$, it allows us to determine whether $A=J(C)$, while given only a ppav $(A,\T)$, this does not provide a way to construct such a curve.

It is then natural to ask what happens for higher-dimensional subvarieties of minimal class. Debarre \cite{debarremin} proved that for any $d$ the loci $W^d(C)$ and $-W^d(C)$ are the only subvarieties of $J(C)$ of the minimal class, and one can ask whether their existence also characterizes Jacobians. It turns out that for intermediate Jacobians of cubic threefolds (these are in $\cA_5$) the Fano surface of lines also has minimal cohomology class (see \cite{clgr,beauvilleprymsubvar}, and also \cite{ranmin2}). Motivated also by works of Ran, Debarre made the following
\begin{conj}[\cite{debarremin}]
If a ppav $(A,\T)$ has a $d$-dimensional subvariety of minimal class, then it is either a Jacobian of a curve or a 5-dimensional intermediate Jacobian of a cubic threefold.
\end{conj}
Debarre proves that this indeed gives a weak solution to the Schottky problem, i.e.~that $\cJ_g$ is an irreducible component of the locus of ppavs for which there exists a subvariety of the minimal cohomology class. While the conjecture remains open, an exciting new approach relating the minimal classes to generic vanishing was recently introduced by Pareschi-Popa \cite{papomin}, who also made a refinement of the above conjecture involving the theta-dual of a subvariety of a ppav.

\begin{rem}
In analogy with the situation for Jacobians, for the Prym-Schottky problem one is naturally led to look at the Abel-Prym curve $\widetilde C\to J(\widetilde C)\to Prym(\widetilde C\to C)$. It can be shown that this has twice the minimal cohomology class, i.e.~$2\frac{[\T]^{g-1}}{(g-1)!}$, and one wonders whether an analog of Matsusaka-Ran criterion holds in this case. Welters \cite{weltersprymcurve} showed that Pryms are an indeed an irreducible component of the locus of ppavs for which there exists a curve representing twice the minimal class, and described the other components of this locus. It is then natural to try to study maps of $Sym^d(\widetilde C)$ to the Prym, and the cohomology classes of their images, but it seems nothing is known here.
\end{rem}
To summarize, this approach to the Schottky problem gives a complete geometric solution to the weaker version of the problem: determining whether a given ppav is the Jacobian of a given curve.

\section{Projective embeddings of a ppav: the geometry of the Kummer variety}\label{kummerimage}
Another approach to the Schottky problem is by embedding a ppav (in fact its quotient under $\pm 1$) into a projective space, and studying the properties of the image. We recall that by definition of a ppav $(A,\T)$ the line bundle $\T$ is ample, but since $h^0(A,\T)=1$, clearly $\T$ is not very ample. The Lefschetz theorem (see \cite{bila}) states that in fact $n\T$ is very ample for any $n\ge 3$, while $2\T$ is a very ample line bundle on $A/\pm 1$. Indeed, recall that on $A_\tau$ a basis of sections of $2\T_\tau$ is given by theta functions of the second order $\lbrace\T[\e](\tau,z)\rbrace$, for all $\e\in\frac12\bZ^g/\bZ^g$, given by formula (\ref{2Tdef}). Recall that all of these are even functions of $z$.
\begin{df}[Kummer map and Kummer variety]
The Kummer map is the embedding
\begin{equation}\label{kum}
Kum:A_\tau/\pm 1\hookrightarrow\bP^{2^g-1};\quad Kum(z):=\lbrace \T[\e](\tau,z)\rbrace_{\rm all\ \e\in\frac12\bZ^g/\bZ^g}
\end{equation}
and we call its image $Kum(A_\tau/\pm1)\subset\bP^{2^g-1}$ the Kummer variety. Notice that the involution $\pm 1$ has $2^{2g}$ fixed points on $A$, which are precisely $A[2]$, and thus the Kummer variety is singular at their images in $\bP^{2^g-1}$.
\end{df}
The Kummer variety is a $g$-dimensional subvariety of $\bP^{2^g-1}$ and one can ask how general it is. The following striking result shows that the Kummer image of a Jacobian is very special (to prove this one uses Riemann's theta singularity theorem to prove Weil reducibility --- that the intersection of $\T$ with some translate is reducible --- and then uses Koszul cohomology and Riemann's bilinear relation to obtain an equivalent statement in terms of the Kummer map: see \cite{acgh,tasurvey} for discussions):
\begin{thm}[Fay(-Gunning)'s trisecant identity \cite{faybook,gunninggen}]
For any $p,p_1,p_2,p_3\in C$ the following three points on the Kummer variety are collinear:
\begin{equation}\label{trisecant}
Kum(p+p_1-p_2-p_3),Kum(p+p_2-p_1-p_3),
Kum(p+p_3-p_1-p_2)
\end{equation}
(where we view these points as line bundles, i.e.~in $\Pic^0(C)=J(C)$)
\end{thm}
This result can be reformulated from the point of view of the geometry of the complete linear series $|2\T|$. Indeed, for any $z_1,z_2,z_3\in A$ consider the vector subspace of $H^0(A,2\T)$ consisting of sections vanishing at these points. Since $Kum$ is the map given by $H^0(A,2\T)$ the number of conditions imposed on a section of $2\T$ by vanishing at these points is equal to the dimension of the linear span of $\lbrace Kum(z_i)\rbrace$. Thus the trisecant identity states that on a Jacobian if the three points are obtained from the points on $C$ as above, the vanishing condition at these 3 points on $J(C)$ fails to impose independent conditions on $|2\T|$.

\smallskip
The trisecant identity implies that the Kummer image of a Jacobian admits a 4-dimensional family of trisecant lines (and thus is very far from being in general position in $\bP^{2^g-1}$). It is natural to wonder whether this is a characteristic property of Kummer images of Jacobians.
\begin{thm}[Gunning \cite{gunningtrisecant}]
If for a ppav $(A,\T)\in\cA_g^{ind}$ there exist three fixed points $p_1,p_2,p_3\in A$ (with some mild general position assumption) such that there exist infinitely many $p$ such that (\ref{trisecant}) is a trisecant of the Kummer variety, then $A\in\cJ_g$ (and moreover there exists a translate of $C\subset J(C)$ such that $p,p_1,p_2,p_3$ are contained in the image).
\end{thm}
The proof of this result proceeds by reduction to the Matsusaka-Ran criterion. Similarly to that result, this theorem gives a complete solution to a weaker version of the Schottky problem: taking the Zariski closure of all $p$ satisfying (\ref{trisecant}) already gives a curve $C\subset A$, and the criterion can be used to verify whether $(A,\T)=J(C)$.

It is thus natural to ask whether one could generalize the statement of the trisecant identity in such a way as to get a characterization of Jacobians not involving a curve contained in the ppav to start with. The general multisecant identity (which is unfortunately less well-known than the trisecant identity) is as follows
\begin{thm}[Gunning \cite{gunningidentities}]
For any curve $C\in\cM_g$, for any $1\le k\le g$ and for any $p_1,\ldots, p_{k+2}, q_1, \ldots, q_k\in C$ the $k+2$ points of the Kummer variety \begin{equation}
Kum\left(2p_j+\sum_{i=1}^k q_i-\sum_{i=1}^{k+2} p_i\right),\qquad j=1\dots k+2
\end{equation}
(where we again identify $J(C)=\Pic^0(C)$) are linearly dependent.
\end{thm}
Similarly to the above, this is the statement that the $k+2$ points on the Jacobian constructed as above fail to impose independent conditions on the linear system $|2\T|$. By fixing $p_1,\ldots,p_{k+2}$ and varying $q_1,\ldots,q_k$ we can further interpret this as the $k+2$ fixed points $$\left\lbrace 2p_j-\sum\limits_{i=1}^{k+2}p_i\right\rbrace_{j=1\ldots k+2}$$ failing to impose independent condition on $|2\T|$ translated by an arbitrary element $q_1+\ldots+q_k\in W^k(C)$.

Now we recall that the map $W^g(C)\to J(C)$ is surjective, and thus the $k=g$ case of the multisecant identity above is the statement that some $g+2$ points on $J(C)$ fail to impose independent conditions on any translate of $|2\T|$. This turns out to also be interpretable in terms of some addition theorems for Baker-Akhiezer functions, and the following result was conjectured by Buchstaber-Krichever \cite{bukr1,bukr2}:
\begin{thm}[Grushevsky \cite{gcubics}, Pareschi-Popa \cite{paposchottky}]
For any $(A,\T)\in\cA_g^{ind}$ and $p_1,\ldots, p_{g+2}\in A$ (in general position, which is a different condition for the two proofs!), if
$$
 \forall z\in A\qquad \lbrace Kum(2p_i+z)\rbrace_{i=1\ldots g+2}\subset\bP^{2^g-1}
$$
are linearly dependent, then $A\in \cJ_g$ (and moreover all $p_i$ lie on a translate of $C\subset J(C)$).
\end{thm}
These results give a weak solution (in the sense of possible extra components that appear if we ignore the general position assumption) to the original Schottky problem --- i.e.~we do not start with a curve $C\subset A$ given. Still the proof proceeds by reduction to Gunning's trisecant criterion above: essentially one considers the loci of all translates of $|2\T|$ on which a subset of less than $g+2$ of the given points fails to impose independent conditions, and once we get to three points not imposing independent conditions, this is a trisecant. Note that this gives a potentially new way of identifying the image of $C$ inside $J(C)$. However, this characterization of Jacobians is still not easy to use explicitly, both because of the general position assumption and because of the necessity to choose/guess appropriate $g+2$ points in a ppav.

Another approach to the Schottky problem based on the trisecant identity is obtained by degenerating the trisecant. Indeed, if we let all the points $p$ and $p_i$ vary, when all $p_i$ come together, the trisecant will degenerate to a flex line, i.e.~to a line tangent to $Kum(A)$ at a point (actually at $p-p_1$) with multiplicity 3. Varying $p$ and $p_1=p_2=p_3$, we thus get a family of such flex lines, i.e.~a family of partial differential equations for the Kummer image, with parameter $p-p_1$. It was shown by Welters \cite{welters} that the existence of such a family of flex lines also characterizes Jacobians --- this is again a solution to the weaker version of the Schottky problem, as we start with a given curve in a ppav.

Expanding these differential equations in Taylor series in $p-p_1$ near $p=p_1$ (so that the flex line goes through zero), we get an infinite sequence of partial differential equations of increasing orders for the Kummer variety at 0, essentially each equation saying that there exists a $k$-jet of a family of flex lines, assuming there exists a $(k-1)$-jet of a family of flexes. This sequence of equations is called the Kadomtsev-Petviashvili (KP) hierarchy. Since we know that for Jacobians there exists a family of flex lines of the Kummer variety, each term of this Taylor expansion, i.e.~each equation of the KP hierarchy must be satisfied by the Kummer image of a Jacobian. We refer to \cite{dusurvey,tasurvey,bukrsurvey} for more on this circle of ideas, and for a much more detailed introduction to the integrable systems point of view on this problem.

Since we are in an algebraic setting, conversely it is natural to expect that if the Kummer variety of a ppav satisfies all the equations of the KP hierarchy, i.e.~if it admits an $\infty$-jet of a family of flex lines, then there exists an actual family of flex lines, and thus the ppav is a Jacobian by the result of Welters. Moreover, from algebraicity one could expect that a finite number of equations from the KP hierarchy would suffice; this was proven by Arbarello-De Concini \cite{ardcequations} to indeed be the case.

Novikov conjectured that in fact the first equation of the KP hierarchy would already suffice to characterize Jacobians, i.e.~that if for some (indecomposable) ppav there exists a 1-jet of a family of flex lines of the Kummer variety at 0, then the ppav is a Jacobian. Explicitly this first equation of the KP hierarchy, called the KP equation, written in terms of theta functions of the second order, has the following form \cite{dubrovinKP}:
\begin{equation}\label{KP}
\left(\frac{\partial^4}{\partial U^4} +\frac{3}{4}\frac{\partial^2}{\partial V^2}-\frac{\partial^2}{\partial U\partial W}+C\right)\T[\e](\tau,z)|_{z=0}=0
\quad \forall \e\in\frac12\bZ^g/\bZ^g,
\end{equation}
for some $U,V,W\in\bC^g,C\in\bC$. Essentially this equation is the condition that the Kummer image of 0 and the second and fourth derivatives of the Kummer image at 0 are linearly dependent --- though the KP is more restrictive, as the coefficients of this linear dependence are not arbitrary. Novikov's conjecture was proven:
\begin{thm}[Shiota \cite{shiota}]
If the Kummer variety of a ppav $A\in\cA_g^{ind}$ admits a 1-jet of a family of flex lines at 0, i.e.~if there exist $U,V,W,C$ such that equation (\ref{KP}) is satisfied, then $A\in\cJ_g$.
\end{thm}
This result is much stronger than the previous theorems for the following two reasons: there is no general position assumption, and to apply this theorem we do not need start with a curve in $A$ (though we still need to choose the parameters $U,V,W,C$ in the KP equation).

One can then wonder whether a still stronger version of this characterization of Jacobians may hold --- whether the existence of just one trisecant (or one flex line) may suffice. Welters \cite{welters} boldly conjectured that the existence of one trisecant of the Kummer variety already guarantees that a ppav is a Jacobian. There have been numerous results in this direction, proving the conjecture under various general position assumptions (see Debarre \cite{detriseclines1,detriseclines2}), and the conjecture was recently proven completely, with no general position assumption:
\begin{thm}[Krichever \cite{kricheverflex,krichevertrisecant}]
For a ppav $A\in\cA_g^{ind}$, if $Kum(A)\subset\bP^{2^g-1}$ has one of the following
\begin{itemize}
\item a trisecant line
\item a line tangent to it at one point, and intersecting it at another (this is a semidegenerate trisecant, when two points of secancy coincide)
\item a flex line (this is a most degenerate trisecant when all three points of secancy coincide)
\end{itemize}
such that none of the points of intersection of this line with the Kummer variety are in $A[2]$ (where $Kum(A)$ is singular), then $A\in\cJ_g$.
\end{thm}
Here we still need to be able to choose the 3 points of trisecancy or the parameters of a line, but there is no curve involved, and for a non-degenerate trisecant we do not even need to consider the derivatives of the Kummer map. It is thus natural to expect that this result would be much harder to prove, as the jet of the curve needs to be constructed from scratch. There are many conceptual and technical difficulties in the proof, and methods of integrable systems are used extensively.

\smallskip
For the Prym-Schottky problem, it is then natural to wonder whether the Kummer varieties of Pryms have any special geometric properties. Indeed, the following result was obtained
\begin{thm}[Quadrisecant identity: Fay \cite{fayprym}, Beauville-Debarre \cite{bedeprym}]
For any $p,p_1,p_2,p_3\in \widetilde C\to Prym(\widetilde C\to C)$ on the Abel-Prym curve the four points of the Kummer variety
\begin{equation}\label{4secant}
\begin{matrix}
&Kum(p+p_1+p_2+p_3),&Kum(p+p_1-p_2-p_3),\\
&Kum(p+p_2-p_1-p_3),&Kum(p+p_3-p_1-p_2)
\end{matrix}
\end{equation}
lie on a 2-plane in $\bP^{2^g-1}$.
\end{thm}
\begin{rem}
Note that there is something slightly puzzling here --- while for the trisecant identity all the points were naturally in $J(C)=\Pic^0(C)$, here the degrees of the four divisors on $\widetilde C$ are $2,0,0,0$, so that one needs to clarify what exactly is meant by the points above. Analytically one can lift the entire discussion to the universal cover $\cH_g\times\bC^g$ of the universal family; we refer to \cite{mumfordprym,bedeprym,casalaina} for more on divisors on Prym varieties.
\end{rem}

Analogously to Gunning's characterization of Jacobians by the existence of a family of trisecants, we have the following
\begin{thm}[Debarre \cite{debarrestratification}]
If for a ppav $(A,\T)\in\cA_g^{ind}$ there exist three fixed points $p_1,p_2,p_3\in A$ (with some mild general position assumption) such that there exist infinitely many $p$ such that (\ref{4secant}) is a quadrisecant 2-plane of the Kummer variety, then $A\in\overline{\cP_g}$ (and moreover $p,p_1,p_2,p_3$ lie on an Abel-Prym curve).
\end{thm}
Following the outline of our discussion for Jacobians, it is now natural to ask whether one could increase the number of points and get a suitable multisecant identity for Pryms. It appears that by generalizing the methods of \cite{fayprym,bedeprym,gunningidentities} one could indeed construct a $(2k+2)$-dimensional family of $(2k+2)$-secant $2k$-planes of the Kummer images of Pryms (for comparison, Gunning's multisecant identity for Kummer varieties of Jacobians gives a $(2k+2)$-dimensional family of $(k+2)$-secant $k$-planes). In particular, to try to emulate the results of \cite{gcubics,paposchottky} one would need to consider $(2g+2)$-secant $2g$-planes --- equivalently sets of $2g+2$ points imposing linearly dependent conditions on an arbitrary translate of $|2\T|$ (for comparison, for Jacobians we had $g+2$ points). This gets very complicated, and in particular the following questions seems intriguing and wide open:
\end{rem}
\begin{OP}
For some $A\in\cA_g^{ind}$ suppose there exist $k$ points $p_1,\ldots,p_k\in A$ in general position imposing less than $k$ independent conditions on any translate of $|2\T|$ (equivalently, such that for any $z\in A$ the Kummer images $\lbrace Kum(p_i+z)\rbrace$ are linearly dependent). For which maximal $k$ does this imply that $A\in\overline{\cP_g}$ (we know that $k=g+2$ implies $A\in\cJ_g\subset\overline{\cP_g}$)? It is likely that the  maximum that $k$ could be is at most $2g+2$, so what loci within $\cA_g^{ind}$ are characterized by this property for $g+2<k<2g+2$? Do they give examples of some natural loci $\cJ_g\subsetneq X\subsetneq\overline{\cP_g}$, or do Pryms satisfy this condition with $k\le g+5$ (recall that for $g\ge 5$ we have $\dim\overline{\cP_g}=3g=3+(3g-3)=3+\dim\cM_g$)?
\end{OP}
Of course the above characterization of Pryms by a family of quadrisecant planes makes one wonder whether one could construct an integrable hierarchy of partial differential equations satisfied by the theta function of Pryms, by allowing the quadrisecant plane to degenerate appropriately. This is indeed possible, and one can thus obtain the Novikov-Veselov, the BKP, and the Landau-Lifschitz hierarchies. One then naturally expects finitely many of the equations in these hierarchies to characterize Pryms --- which is indeed the case. However, the analog of the Novikov conjecture --- that the first equation of one of these hierarchies characterizes Pryms --- is open, see \cite{taimanovdiffeq}. It is only known (Shiota \cite{shiotaprym,shiotaprym2}) that the BKP equation characterizes Pryms under a certain general position assumption. We refer to \cite{tasurvey} for a detailed survey of what is known about the Prym-Schottky problem and soliton equations.

\smallskip
Similarly to the trisecant conjecture for Jacobians, one can wonder whether the existence of one quadrisecant 2-plane of the Kummer variety may characterize Pryms. This turns out not to be the case: Beauville-Debarre \cite{bedeprym} construct an example of $A\in\cA_g^{ind}\setminus\overline{\cP_g}$ such that $Kum(A)$ has a quadrisecant 2-plane. Thus no analog of quadrisecant conjecture has been proposed for Pryms, though Beauville-Debarre proved in \cite{bedeprym} that $\cP_g$ is an irreducible component of the locus in $\cA_g$ consisting of ppavs whose Kummer varieties admit a quadrisecant 2-plane.

A suitable analog of the trisecant conjecture was recently found for Pryms using ideas of integrable systems; it admits a clear geometric formulation:
\begin{thm}[Grushevsky-Krichever \cite{grkr}]
If for some $A\in\cA_g^{ind}$ and some $p,p_1,p_2,p_3\in A$ the quadrisecant condition (\ref{4secant}) holds, and moreover there exists another quadrisecant given by (\ref{4secant}) with $p$ replaced by $-p$, then $A\in\overline{\cP_g}$.
\end{thm}
One can think of this pair of quadrisecants in the following way: recall that $Prym(\widetilde C\to C)$ is the connected component of zero in the kernel of the map $J(\widetilde C)\to J(C)$, and that this kernel has two connected components. The theorem essentially says that both the Prym and the other component of this kernel must admit a quadrisecant.

We note that similarly to the formulation of the trisecant conjecture, this theorem is a complete (not weak) solution to the Prym-Schottky problem. There is no general position assumption here (which would result in extra components), and one does not start with the Abel-Prym curve. The proof of this result is very involved and again uses the ideas coming from integrable systems to construct first an infinitesimal jet of $\widetilde C\to C$ and then argue that there indeed exists such a cover.

\begin{rem}
One possible very interesting question to which one could try to apply the above theorem is to try to approach the Torelli problem for Pryms: indeed, the map from the moduli space of \'etale double covers of curves of genus $g$ to $\cA_{g-1}$ is known to be generically injective for $g\ge 7$ (see \cite{dosm} for more on the Prym map), but it is never injective, since Donagi's tetragonal construction \cite{donagitetragonal} can be used to construct different covers with isomorphic Pryms. It was conjectured (see Donagi \cite{donagitetragonal}, Debarre \cite{debarreprymtorelli}) that for $g\ge 11$ this was the only cause of non-injectivity of the Prym Torelli map. However, a recent preprint of Izadi and Lange \cite{izadilange} constructs examples of curves of arbitrary Clifford index (thus the tetragonal construction does not apply) such that the Prym Torelli map is not injective at the corresponding Prym. It is tempting to try to apply the above characterization of Pryms to study this problem: the data of the two quadrisecants in fact recovers the cover $\widetilde C\to C$, not only the Prym, and thus the question is whether a given Prym may have more than one family of quadrisecants.
\end{rem}

\smallskip
To summarize, the approach to the Schottky problem via the geometry of the Kummer variety yields complete (strong --- no extra components) solutions to the Schottky and Prym-Schottky problems. These solutions do not require an a priori knowledge of a curve, as does the minimal cohomology class approach. However, there is still a finite amount of data involved --- choosing the points of secancy. Thus given an explicit ppav it is not yet clear how to apply these solutions of the Schottky problem or to relate this approach to the Schottky-Jung approach; in particular it does not yet seem possible to apply the results in this section to solve either open problem \ref{SchottkyOP} or open problem \ref{ColemanOP}.

\section{$\Gamma_{00}$ conjecture}\label{gamma00}
In the previous section we saw that the geometry of the Kummer variety, i.e.~the geometry of the $|2\T|$ linear system, is of great importance in studying the Schottky problem. The second order theta functions $\T[\e](\tau,z)$ seem to be central to the study of the Schottky problem: their values at $z=0$, for $\tau$ varying, are theta constants giving the embedding $Th:\cA_g(2,4)\to\bP^{2^g-1}$, which is the subject of the classical and the Schottky-Jung approaches (see sections \ref{classical},\ref{schottkyjung}); while for $\tau$ fixed and $z$ varying theta functions of the second order define the map $Kum:A_\tau/\pm 1\hookrightarrow\bP^{2^g-1}$, properties of which we discussed in section \ref{kummerimage} (in fact the intersection of the images of $Th$ and $Kum$ in $\bP^{2^g-1}$ is also of interest --- see \cite{vgvdg}).

The difficulty in trying to relate the geometric properties of the Kummer image (eg. the trisecants), or other properties of the theta divisor --- eg. the singularities, as in section \ref{andreottimayer}, and algebraic equation for theta constants is that it is not clear how to capture the properties of $\T$ and $Kum(A)$ just from $Th(A_\tau)=Kum(0)$.

\begin{rem}
An interesting problem showing that we are yet unable to translate geometric conditions for the theta divisor into modular forms is as follows. Casalaina-Martin and Friedman \cite{cmfr,casalaina2} showed that intermediate Jacobians of cubic threefolds (which we already mentioned in relation to the Schottky-Jung approach) are characterized within $\cA_5$ by the property of having a single point of multiplicity 3 on $\T$ (this point is thus in $A[2]^{odd})$. However, the following question remains open; solving it could also yield a solution to open problem \ref{SchottkyOP}, as it would potentially allow one to further refine the Schottky-Jung approach in genus 5, dealing with the big Schottky-Jung locus, and in particular possibly proving conjecture \ref{sjconj} for $g=5$:
\begin{OP}
Determine the defining ideal for $\overline{Th(\cI(2,4))}\subset \overline{Th(\cA_5(2,4))}$, where $\cI\subset\cA_5$ denotes the locus of intermediate Jacobians of cubic threefolds.
\end{OP}
We note that it is also very interesting to try to understand the geometry of the boundary of $\cI$, and work on this has been done by Casalaina-Martin and Laza \cite{cmla}; see also \cite{grsmordertwo}.
\end{rem}

\smallskip
Thus it is natural to try to study the properties of $|2\T|$ at $0\in A$. Note that the statement of the trisecant conjecture explicitly excludes $0$, because the Kummer variety is singular at $Kum(0)$. One thus wonders what an appropriate notion of a degenerate trisecant through $Kum(0)$ might be, and it is natural to expect that it would involve the tangent cone to the Kummer variety at 0.
\begin{df}[van Geemen-van der Geer \cite{vgvdg}]
The linear system $\Gamma_{00}\subset |2\T|$ is defined to consist of those sections that vanish to order at least 4 at the origin:
$$\Gamma_{00}=\lbrace f\in H^0(A,2\T) \mid mult_0f\ge 4\rbrace.$$
\end{df}

We note that since all sections of $2\T$ are even functions of $z$, all first and third partial derivatives vanish automatically, and thus the condition for $f$ to lie in $\Gamma_{00}$ is really for the $f$ and its second partial derivatives to vanish at zero. We recall now that theta functions of the second order generate $H^0(A,2\T)$, so that
$$
f(z)=\sum\limits_{\e\in\frac12\bZ^g/\bZ^g} c_\e\T[\e](\tau,z),
$$
and thus the condition $f\in\Gamma_{00}$ is equivalent to the following system of $\frac{g(g+1)}{2}+1$ linear equations for the coefficients $\lbrace c_\e\rbrace$:
$$
 0=\sum c_\e\T[\e](\tau,0)=\sum c_\e\frac{\partial^2\T[\e](\tau,z)}{\partial z_a \partial z_b}|_{z=0}\qquad\forall 1\le a\le b\le g.
$$
One can show (eg. \cite{dubrovinKP}) that for any ppav in $\cA_g^{ind}$ the rank of this linear system of equations is maximal, and thus that $\dim \Gamma_{00}=2^g-\frac{g(g+1)}2-1$ on any indecomposable ppav (in particular it is empty for $g=1,2$, and one-dimensional for $g=3$).

\begin{rem}
The condition $f\in\Gamma_{00}$ can actually be reformulated purely in terms of $\tau$, not involving $z$. The Riemann theta function satisfies the heat equation, and so do the theta functions of the second order, for which we have
$$
 \frac{\partial^2}{\partial z_i \partial z_j}\T[\e](\tau,z)= \frac{2\pi i}{1+\de_{a,b}}\frac{\partial}{\partial \tau_{ab}}\T[\e](\tau,z)
$$
(this can be easily verified from definition, by differentiating the Fourier series term by term). Thus one can study the linear system $\Gamma_{00}$ purely in terms of theta constants of the second order and their derivatives (with respect to $\tau$), which may give a hope of eventually relating this to modular forms and the Schottky-Jung approach.
\end{rem}

A priori the base locus $Bs(\Gamma_{00})$ contains 0, and it is not clear why it would contain other points. From the above description of $\Gamma_{00}$ as a system of linear equations on coefficients $c_\e$ it follows that the condition $z\in Bs(\Gamma_{00})$ is equivalent to
$$
Kum(z)\in \left\langle Kum(0), \frac{\partial}{\partial\tau_{ab}} Kum(0)\right\rangle_{\rm linear\ span},
$$
i.e.~to the existence of $c,c_{ab}\in\bC$ such that
\begin{equation}\label{dependence}
  \T[\e](\tau,z)=c \T[\e](\tau,0)+\sum\limits_{1\le a\le b\le g} c_{ab}\frac{\partial\T[\e](\tau,0)}{\partial\tau_{ab}}\quad \forall\e\in\frac12\bZ^g/\bZ^g.
\end{equation}
Indeed, $z\in Bs(\Gamma_{00})$ means that $z$ imposes no conditions on $\Gamma_{00}$, and thus that the vanishing of $f(z)$ must be a consequence of $f$ satisfying the linear conditions defining $\Gamma_{00}$.

We note that this condition is indeed very similar to the condition of the existence of a semidegenerate trisecant of a Kummer variety --- in essence it is saying that a line through $Kum(z)$ and $Kum(0)$ lies in the tangent cone to the Kummer variety at 0. However, we are unaware of any formal relationship or implication between this statement and that of a semidegenerate trisecant through a point not of order two.

For Jacobians one can use Riemann's theta singularity theorem combined with Riemann's bilinear addition formula to show that the difference variety $C-C:=\lbrace p-q|p,q\in C\rbrace\subset \Pic^0(C)=J(C)$ is contained in $Bs(\Gamma_{00})$ --- this is essentially due to Fr\"obenius \cite{frobenius}, see also \cite{vgvdg}. Moreover, it turns out that for $p,q\in C\subset J(C)$ the rank of the matrix of coefficients $c_{ab}$ in (\ref{dependence}) is equal to one.

It turns out that this is the entire base locus:
\begin{thm}[Welters \cite{weltersC-C} --- set-theoretically, Izadi \cite{izadiC-C} --- scheme-theoretically]
For any $g\ge 5$ and any $C\in\cM_g$ we have on $J(C)$ the equality $Bs(\Gamma_{00})=C-C$.
\end{thm}

One then wonders whether this property characterizes Jacobians:
\begin{conj}[$\Gamma_{00}$ conjecture; van Geemen-van der Geer \cite{vgvdg}]
If $Bs(\Gamma_{00})\ne\lbrace0\rbrace$ for some $(A,\T)\in\cA_g^{ind}$, then $A\in\cJ_g$.
\end{conj}
This conjecture was proven for $g=4$ by Izadi \cite{izadiA4}, for a generic Prym for $g\ge 8$ by Izadi \cite{izadiprymgamma00}, and for a generic ppav for $g=5$ or $g\ge 14$ by Beauville-Debarre-Donagi-van der Geer \cite{bddvdg}. However, in order to obtain a solution of the Schottky problem, proving the $\Gamma_{00}$ conjecture for a generic ppav is insufficient --- a generic ppav is not a Jacobian anyway. At the moment there does not seem to be a promising algebro-geometric approach to the $\Gamma_{00}$ conjecture in general, but the following result was recently obtained by using integrable systems methods:
\begin{thm}[Grushevsky \cite{grgamma00}]
If for some $A\in\cA_g^{ind}$ the linear dependence (\ref{dependence}) holds with $\operatorname{rk}(c_{ab})=1$, then $A\in\cJ_g$.
\end{thm}
The proof of this result actually exhibits the similarity between the $\Gamma_{00}$ conjecture and Krichever's characterization of Jacobians by a semidegenerate trisecant --- both of these conditions turn out to imply the same differential-difference equation for the theta function (which is what eventually characterizes Jacobians analytically).

\begin{rem}
The $\Gamma_{00}$ conjecture turns out to also be related to the Schottky-Jung theory. Indeed, one boundary component of the Deligne-Mumford compactification $\overline{\cM_g}$ is the moduli space of pointed curves $\cM_{g-1,2}$. Given a modular form vanishing identically on $\overline{\cJ_g}$, one can consider its restriction to the boundary. After a suitable blowup it means that one is considering the next term of the Taylor expansion (called Fourier-Jacobi expansion in this case) of this modular form near the boundary of $\overline{\cM_g}$, and the vanishing of this modular form would imply --- and may be implied --- by the identical vanishing of this next term in the expansion along $\cM_{g-1,2}$, which is the question of a certain modular function $F(\tau,z)$ lying in $\Gamma_{00}$ for all $\tau$ (see \cite{grsmtwopoint,grsmgen5} for a demonstration of how this works, and applications to some questions in perturbative string theory). Mu\~noz-Porras suggested using $\Gamma_{00}$ conjecture as an approach to proving conjecture \ref{sjconj} by such degeneration methods, but this has not yet been accomplished.
\end{rem}

The $\Gamma_{00}$ linear system also turns out to be conjecturally related to another approach to the Schottky problem. This line of investigation started probably from the work of Buser-Sarnak \cite{busa} who showed that for $g$ large, the length of the shortest period of a Jacobian (i.e.~the vector in the lattice $\bZ^g+\tau\bZ^g$ with minimal length) is much less than for a generic ppav. This immediately gives an effective combinatorial way to show that some ppavs are not Jacobians --- if their periods are too long --- but being an open condition of course cannot characterize Jacobians. In this direction, of finding explicit conditions guaranteeing that a ppav is {\it not} a Jacobian, there has also been recent progress in working with ppavs admitting extra automorphisms --- note in particular Zarhin's work \cite{zarhin1,zarhin2} where examples of ppavs with automorphisms that are not Jacobians are constructed.

The length of the shortest period was related to the value of the Seshadri constant by Lazarsfeld \cite{lazarsfeld}, and Seshadri constants on abelian varieties were further studied
by Nakamaye \cite{nakamaye}, Bauer \cite{bauer}, and Bauer-Szemberg \cite{basz}. It turns out that for $g$ large the Seshadri constant for any Jacobian is much smaller than for a generic ppav of that dimension. Unlike the length of the periods, which cannot possibly characterize Jacobians, it turns out that the Seshadri constant may actually give a solution at least to the hyperelliptic Schottky problem, and is related to $\Gamma_{00}$:
\begin{thm}[Debarre \cite{debarreseshadri}]
If the $\Gamma_{00}$ conjecture holds, hyperelliptic Jacobians are characterized in $\cA_g$ by the value of their Seshadri constant.
\end{thm}

\smallskip
It thus seems that the $\Gamma_{00}$ conjecture approach, though not yet culminating even in a weak solution to the Schottky problem, is related to many other approaches to the Schottky problem, and may possibly serve as a tool to further study the geometric properties of the Jacobian locus $\cJ_g\subset\cA_g$ by degeneration methods.

%\bibliographystyle{alpha}
%\bibliography{../../sam_biblio}

\end{document}